\documentclass[leqno,12pt]{article}

\usepackage{amsmath}
\usepackage{amscd}
\usepackage{amsopn}
\usepackage{amsthm}
\usepackage{amsfonts,amssymb}
\usepackage{amsfonts,bbm}
\usepackage{latexsym}

\makeatletter
\@addtoreset{equation}{section}
\makeatother

\newtheorem{lemma}{Lemma}[section]
\newtheorem{proposition}[lemma]{PROPOSITION}
\newtheorem{corollary}[lemma]{COROLLARY}
\newtheorem{theorem}[lemma]{THEOREM}
\newtheorem{remark}[lemma]{REMARK}

\newcommand{\real}{\mathbbm{R}}
\newcommand{\nat}{\mathbbm{N}}

\newcommand{\ganz}{\mathbbm{Z}}
\newcommand{\ganzd}{\mathbbm Z^d}

\renewcommand{\a}{\alpha}
\renewcommand{\b}{\beta}

\newcommand{\ve}{\varepsilon}

\newcommand{\reald}{{\real^d}}

\newcommand{\on}{\quad\text{ on }}

\newcommand{\inv}{^{-1}}
\newcommand{\ov}{\overline}

\newcommand{\M}{\mathcal M}

\newcommand{\Z}{\mathcal Z}
\newcommand{\K}{\mathcal K}
\newcommand{\Q}{\mathcal Q}

\newcommand{\dist}{\mbox{\rm dist}}

\newcommand{\itemframe}%
{\setlength{\parskip}{10pt}\begin{enumerate} \setlength{\topsep}{10pt}%
\setlength{\itemsep}{15pt}\setlength{\parsep}{5pt}}

\newcommand{\vy}{\ve_y}
\newcommand{\vx}{\ve_x}
\newcommand{\nus}{\nu^{A_s}}
\newcommand{\nut}{\nu^{A_t}}

\newcommand{\rhat}{\hat R}

\newcommand{\Px}{\mathcal P(X)}
\newcommand{\Mp}{\M(\Px)}

\newcommand{\mne}{\bigl(\M_\nu(\Px)\bigr)_e}
\newcommand{\mnp}{\M_\nu(\Px)}

\newcommand{\mns}{\M_\nu(S(W))}
\newcommand{\mnh}{\M_\nu(H(W))}
\newcommand{\mxs}{\M_x(S(W))}

\newcommand{\mnse}{ \bigl(\mns\bigr)_e}
\newcommand{\mxse}{ \bigl(\mxs\bigr)_e}

\newcommand{\onu}{\overset{\,\circ}\nu{}}

\newcommand{\wtnu}{{}^W\!\tilde\nu}

\newcommand{\uc}{{U^c}}
\newcommand{\unc}{{U_n^c}}
\newcommand{\vc}{{V^c}}
\newcommand{\wc}{{W^c}}

\renewcommand{\L}{\Lambda_k}
\newcommand{\lan}{\lambda_n}
\newcommand{\sulan}{\sum_{n=1}^k\lan}
\newcommand{\nun}{\nu^{K_n}}
\newcommand{\nuun}{\nu^{U_n\cup\wc}}
\newcommand{\nukn}{\nu^{K_n\cup\wc}}
\newcommand{\nutkn}{\tilde\nu^{K_n\cup\wc}}
\newcommand{\mut}{\tilde\mu}

\title{Convexity of limits of harmonic measures}
\author{WOLFHARD HANSEN and IVAN NETUKA
\thanks{The work is a part of the research project MSM
0021620839 financed by  MSMT.}
}

\date{Preliminary version}

\begin{document}
 
\maketitle
\begin{abstract}
It is shown that, given  a point $x\in\reald$, $d\ge 2$, 
and open sets $U_1,\dots,U_k$ in~$\reald$ containing $x$, 
any convex combination 
of the  harmonic measures  $\vx^\unc$ for~$x$ with respect to $U_n$,
$1\le n\le k$, is the limit of a sequence $(\vx^{W_m^c})_{m\in\nat}$ 
of~harmonic measures,
where each $W_m$ is an open subset of $U_1\cup\dots \cup U_k$
containing~$x$. This answers a question raised in connection 
with Jensen measures.

 More 
generally, we prove that, for arbitrary measures on an open set $W$,
the set of extremal representing measures, with respect to 
the cone of continuous potentials on $W$ or with respect to the
cone of continuous functions on $\ov W$ which are superharmonic $W$, 
is dense in the compact convex set of all representing measures.

This is achieved approximating balayage on open sets
by balayage on unions 
of balls which are pairwise disjoint 
and very small with respect to their mutual distances and then
shrinking  these balls in a suitable manner.

The results are presented simultaneously for the classical case 
and for the theory of Riesz potentials. 

Finally, a characterization of all Jensen measures and of all extremal 
Jensen measures is given.

\medskip


\end{abstract}

\section{Introduction and main results}

The principal motivation for this paper is the following 
natural question from classical potential theory which has
been raised explicitly in \cite{cole-ransford} in connection
with Jensen measures: 
Given an open subset $\Omega$ in~$\reald$, $d\ge 2$, and a point
$x\in\Omega$, 
is~the  set of limits  of harmonic measures $\vx^\uc$ 
convex (where the sets $U$ are supposed to be relatively compact open
neighborhoods of $x$ in $\Omega$)?
We shall see that the answer is~\textquotedblleft yes\textquotedblright.
In~fact, we shall prove that even for general measures $\nu$
instead of~$\ve_x$
the extremal representing measures are dense in the compact convex
set of all representing measures  (see Corollary~\ref{main2-int}
and Corollary~\ref{main4-int}). 

Our method of sweeping on families of  disjoint balls, which are very
small with respect to their mutual distances, works as well
for the theory of 
Riesz potentials related to the fractional Laplacian 
$-(-\Delta)^{\a/2}$ on $\reald$, $0<\a<\min\{2,d\}$.
Therefore we shall also cover the case of Riesz potentials from  
the very beginning. We recall that classical potential 
theory of the Laplacian is the limiting case $\a=2$. 
The reader, who is interested in the classical case
only, may neglect this generality. He will hardly notice any difference 
in the presentation except from the additional discussion of the 
``Poisson kernel'' for a ball with respect to Riesz potentials 
(which has  a density with respect to Lebesgue measure
on the complement of the ball).

So we shall deal simultaneously with the 
following two  situations: 
\begin{itemize}
\item Classical case: $\a=2$, $X$ is a non-empty open set in $\reald$, 
$d\ge 2$,  such that $\reald\setminus X$ is non-polar, if $d=2$.
\item Riesz potentials: $\a<2$, $X$ is a non-empty open set in
  $\reald$, $d\ge 1$, $d>\a$.
\end{itemize}
Given  $Y\subset X$, $Y^c:=X\setminus Y$ will denote the complement
of $Y$ \emph{with respect to} $X$.

Let $\M(X)$ be the set of all (Radon) measures on $X$, let $\mathcal P(X)$
 denote the convex cone of all continuous real potentials on $X$,
and   $\Mp$ be the set of all $\nu\in\M(X)$ such that $\mu(p)<\infty$
for \emph{some} strictly positive   $p\in \Px$. Let us note that every
finite  measure  on $X$ and hence every $\nu\in\M(X)$ with
compact support is  contained in $\Mp$.
For every $\nu\in\Mp$ and for every  subset $A$ of $X$, let
 $\nu^A$ denote the balayage of  $\nu$ on $A$ \emph{with respect to} $X$,
that is, for every superharmonic function $u\ge 0$ on $X$,
$$
    \nu^A(u):= \int u \,d\nu^A=\int \rhat_u^A\,d\nu,
$$
where $\rhat_u^A(x)=\liminf_{y\to x} R_u^A(y)$ and $R_u^A$ is the infimum
of all positive superharmonic functions $v$ on $X$ majorizing $u$ on $A$.
In particular, 
\begin{equation}\label{bal-est}
              \nu^A(u)\le \nu(u)
\end{equation} 
for every superharmonic function $u\ge 0$ on $X$.
For every $x\in X$, let $\vx$ denote the Dirac measure at $x$. 
It is easily seen that $\nu^A=\int \vx^A\,d\nu(x)$. If $A$ is closed
and $x\in A^c$, then  $\vx^A$ is the restriction of the harmonic 
measure for $x$ and the open set $X\setminus A$ on~$A$. 
Given $A\subset X$, there exists a Borel set  (even a
$G_\delta$-set) $\tilde A$ containing $A$ such that $\nu^{\tilde A}=
\nu^A$ for every $\nu\in\Mp$. To discuss extremal representing
measures for $\nu\in\Mp$ we shall also need reduced measures
$\onu^A$ for Borel sets $A\subset X$. They are defined by
$$
            \int u\,d\onu^A=\int R_u^A\,d\nu, 
$$
$u\ge 0$ superharmonic on  $X$, 
and related to $\nu^A$ by $\onu^A=\nu|_A+(\nu|_{A^c})^A$, since
$R_p^A=p$ on $A$ and $R_p^A=\rhat_p^A$ on $A^c$. If $A$ is open or,
more generally, if $A$ is not thin at any of its points, then
$\onu^A=\nu^A$ for every $\nu\in\Mp$. We refer to \cite{BH} for details.

Let $\K(X)$ denote the linear space of all continuous real functions 
on $X$ with compact support. We recall that a sequence $(\mu_m)$
of Radon measures converges \emph{weakly} to a Radon measure $\mu$
on $X$ if $\lim_{m\to\infty} \mu_m(f)=\mu(f)$ for every $f\in \K(X)$. 
It~is this
convergence for Radon measures we shall use.

Let us fix a natural number $k\ge 2$ and define
$$
            \L:=\{\lambda\in[0,1]^k\colon \sum_{n=1}^k \lambda_n=1\}.
$$
Our  fundamental result is the following.

\begin{theorem}\label{main-int}
Let  $\nu$ be a measure in $\Mp$  which is supported by an open subset $W$
of $X$, let $U_1,\dots, U_k$ be open subsets of $W$, and $\lambda\in\L$. 
Then there exist  finite unions~$C_m$, $m\in\nat$,  of pairwise disjoint
closed balls in $U_1\cup\dots\cup U_k$ such that 
$$
\lim_{m\to\infty}\nu^{ C_m\cup \wc}\ = \ \sulan \nu^{U_n\cup\wc}.
$$
\end{theorem}

The key to Theorem \ref{main-int} is the following result concerning balayage
on finite families of small balls where, given $\gamma\in[0,1]$ and a
closed ball $B$  with center $x$ and radius $r$, the ball  with 
center~$x$ and radius~$\gamma r$ is denoted by $B^\gamma$ 
(see Proposition \ref{essential-2} for a precise formulation).
It will be applied to balayage on subsets $A$ of $W$
\emph{with respect to $W$} in place
of $X$ to deal with balayage measures of the
form $\nu^{A\cup\wc}$.

\begin{proposition}\label{essential-int}
Let $\delta>0$ be small, let $A$ be a union of finitely many
pairwise disjoint closed balls  $B_1,\dots, B_m$ in $X$  which
are sufficiently small with respect to their mutual distances
and to the distance from $\reald\setminus X$, and let $\nu\in\Mp$ 
such that $\nu(A)=0$.
Moreover, let $\lambda\in\L$, let $I_1,\dots,I_k$ be a partition of 
$\{1,\dots,m\}$, 
and $K_n$ be the union of the balls $B_i$, $i\in I_n$, $1\le n\le k$.

Then there exist $\gamma_1,\dots,\gamma_m \in [0,1]$ such that 
$C:=B_1^{\gamma_1}\cup\dots \cup B_m^{\gamma_m}$ satisfies
$$
   \nu^C(B_i)= (1-\delta) \sulan \nu^{K_n}(B_i)
        \qquad \mbox{ for every } 1\le i\le m.
$$
\end{proposition}

  Given two measures $\mu,\nu$ on $X$, we shall write $\mu\prec\nu$
 provided $\mu(p)\le \nu(p)$ for every $p\in \mathcal P(X)$.
For every $\nu\in\Mp$, let $\mnp$ be the set of all measures~$\mu$
on~$X$ such that $\mu(p)\le \nu(p)$ for every $p\in\mathcal P$, that is,
$$
           \mnp=\{\mu\in\Mp \colon \mu\prec \nu\}.
$$
$\mnp$ is a compact convex set with respect to 
weak convergence and its set of extreme points is given by
\begin{equation}\label{ext}
                   \mne=\{\onu^A\colon A \mbox{ Borel subset of } X\}
\end{equation} 
(see \cite{moko-ext} and \cite[VI.12.5]{BH}). Moreover, the 
subset of all $\nu^U$, $U$ open in $X$, as well
as the subset of all $\onu^C$, $C$ compact subset of $X$, 
is dense in $\mne$ (see \cite[VI.1.9]{BH}). Therefore, by the theorem
of Krein-Milman and taking $W=X$, Theorem \ref{main-int} yields
the following.

\begin{corollary}\label{main2-int}
For every $\nu\in\Mp$, $\mne$ is dense in $\mnp$.
\end{corollary} 
 
\begin{remark} 
Let us note that Corollary \ref{main2-int} has the following 
consequence related to Skorokhod stopping 
{\rm (see \cite{rost,falkner,fitzsimmons-sko,falkner-fitzsimmons,
baxter-chacon})}. 
Let $\nu$ be a probability measure on $X$
and let $(X(t))$ be   Brownian motion or an 
$\a$-stable process on $X$ with initial distribution~$\nu$.
Then, for every measure $\mu\prec\nu$,
there exists a sequence $(T_m)$
 of hitting times at relatively compact open subsets 
 $U_m$ of $X$ such that the
 distributions $P^\nu_{X(T_m)}$ converge weakly to 
$\mu$ as $m\to\infty$.
\end{remark}

Next let us consider open subsets  $U_1,\dots,U_k$ of $X$ and 
let $\nu$ be a measure in $\Mp$ which is supported by 
$W:=U_1\cup\dots \cup U_k$.  
Then, for every $1\le n\le k$, there exist  open $(1/m)$-neighborhoods
$\tilde U_{nm}$  of $U_n^c$ in $X$ such that $(\nu^{\tilde U_{nm}})$ 
converges weakly to $\onu^{U_n^c}$ as $m\to\infty$ (see \cite[VI.1.9]{BH}). 
For all $n$ and $m$,  
$\tilde U_{nm}= (\tilde U_{nm}\cap W)\cup \wc$ and 
 $\tilde U_{nm}\cap W$ is an open subset of $W$.
Therefore Theorem \ref{main-int} implies as well the following.

\begin{corollary}\label{main3-int}
Let $U_1,\dots,U_k$ be  open subsets of $X$, $\nu$ be a measure in $\Mp$ 
which is supported by $W:=U_1\cup\dots \cup U_k$, and $\lambda\in \L$. Then
there exist  finite unions $C_m$  of pairwise disjoint closed balls in a
$(1/m)$-neighborhood of $W\setminus (U_1\cap\dots\cap U_k)$ in $W$ 
such that 
\begin{equation}\label{1.4}
\lim_{m\to\infty} \onu^{\wc\cup C_m}\ = \ \sulan \onu^{U_n^c}.
\footnote{In the classical case and for $\nu$ supported by $U_1\cap\dots
\cap U_k$, 
we may choose $C_m$ in a $(1/m)$-neighborhood of $W\cap(\partial
U_1\cup\dots\cup\partial U_k)$.}
\end{equation} 
If $\nu$ is supported by $U_1\cap\dots\cap U_k$, then the reduced
measures may be replaced by balayage measures.
\end{corollary}

Given an open subset $W$ of $X$,
let $S(W)$, $H(W)$ denote
 the set of all continuous  functions on $X$ which are 
$\Px$-bounded (that is, bounded in modulus by  some $p\in\Px$) and
superharmonic  on $W$, harmonic on $W$, respectively.
As for $\Px$, we have sets of representing measures $\mns$ and $\mnh$.
Since semipolar sets are polar and points are polar 
for both the classical case and for Riesz potentials, 
we~see from \cite[VI.9.5]{BH} that the following holds
for Dirac measures $\nu=\vx$ (as customary, we write $\M_x$ instead of 
$\M_{\vx}$);  the proof for the general case $\nu\in \Mp$
will be given in an Appendix.

\begin{theorem}\label{msg}
For every open subset $W$ of $X$ and for every $\nu\in\Mp$ which is supported
by $W$,
\begin{eqnarray*}
\mns&=&\mnp\cap\mnh\\&=&\{\mu\in\mnp\colon \mu^\wc=\nu^\wc\}
                  =\{\mu\in\Mp\colon \nu^\wc\prec\mu\prec \nu\}.
\end{eqnarray*} 
Moreover, $\mns$ is a closed face of $\mnp$ and
$$
     \mnse= \{\onu^A\colon \wc\subset A\subset X,\ A \mbox{ Borel set}\}.
$$
In particular, for every $x\in W$,
$$
     \mxse=\{\vx^A\colon \wc\subset A\subset X\}.
$$
\end{theorem} 

Let us note that, taking $W=X$, we have $\wc=\emptyset$,
$H(W)=\{0\}$, and $S(W)$ is $\Px$. 
Since the measures $\nu^{U\cup\wc}$, $U$ open subset of $W$, are
dense in $\mnse$, Theorem \ref{main-int} also yields the
following.

\begin{corollary}\label{main4-int}
Let $W$ be an open subset of $X$ and let $\nu$ be a measure in $\Mp$  
which is supported by $W$.   Then $\mnse$ is dense in $\mns$.
\end{corollary} 

Finally, restricting our attention to classical potential theory,
let us see how Jensen measures, introduced in function theory in 
\cite{baxter-chacon} and extensively studied in \cite{cole-ransford}
and \cite{ransford}, fit into our considerations.
To~that end we fix an open subset~$\Omega$  of~$\reald$, $d\ge 2$,
 and a point $x\in \Omega$.
A Jensen measure for $x$ with respect to $\Omega$
is probability measure $\mu$ supported on a compact subset of~$\Omega$
such that 
$
               \int u\,d\mu\le u(x)
$
for every superharmonic function $u$ on $\Omega$. Equivalently, since
constants are harmonic and superharmonic functions are increasing
limits of continuous superharmonic functions, the
set $J_x(\Omega)$ of all Jensen measures for $x$ with respect to
$\Omega$ is the set of all Radon measures with compact support
 in~$\Omega$ such that 
$
 \int u\,d\mu\le u(x) 
$
for every continuous superharmonic function~$u$ on $\Omega$.
 Clearly,
Corollary \ref{main3-int} 
implies by \cite[p.\,32]{cole-ransford} that, for every 
$x\in\Omega$, the set of all harmonic measures $\vx^\uc$, 
$U$ open, $x\in U$, $\ov U$~compact in~$\Omega$, is dense in $J_x(\Omega)$
with respect to the weak$^\ast$-topology on $\mathcal C(\Omega)^\ast$,
that is, Question 1.6 in \cite{cole-ransford} has a positive answer.

Further, we shall give the following characterization for Jensen 
measures and extremal Jensen measures.

\begin{theorem}\label{jensen-char}
Let $\Omega$ be an open subset of $\reald$, $d\ge 2$, and $x\in
\Omega$. Then 
$$
  J_x(\Omega)=\bigcup\bigl\{ \mxs\colon W\mbox{ open},\  x\in W,\ 
                             \ov W\mbox{ compact in } \Omega\bigr\}
$$
and
\begin{eqnarray*}
  \bigl(J_x(\Omega)\bigr)_e & =&
\bigcup\bigl\{ \mxse\colon W\mbox{ open},\ x\in W,\
                             \ov W\mbox{ compact in } \Omega\bigr\}\\[2mm]
&=&\bigl\{\vx^{A^c}\colon \ov A \mbox{ compact in }\Omega\bigr\}.
\end{eqnarray*} 
\end{theorem} 

The paper is organized as follows. In the next section we shall approximate
balayage replacing  open sets $U$ by unions of small balls contained
in $U$. In Section 3
we shall prove Proposition \ref{essential-int}.
Section 4 will consist of the proof for Theorem \ref{main-int},
and in Section 5 we shall establish the results on Jensen measures.
The paper is finished by an Appendix where Theorem \ref{msg} is
proven for  general $\nu\in\Mp$.

\section{Approximation of \mathversion{bold}$\nu^{U\cup\wc}$}

Balayage on open sets can be approximated
by balayage on subsets consisting of finitely many balls having 
radii which are \emph{arbitrarily}
 small with respect to their mutual distances (see Proposition \ref{AmU}).
Since this does not seem to be widely known, we include a complete proof.

For every $x\in\reald$ and $r\ge 0$, let $B(x,r)$ denote the
\emph{closed} ball having center $x$ and radius $r$.
Given $x_0\in\reald$,  $a\in (0,1)$, and $m\in\nat$, let 
$$
\Z_m(x_0,a):=
\bigl\{B(z,\frac am)\colon z\in \frac 1m (x_0+\ganz^d)\bigr\},\qquad
Z_m(x_0,a):=\bigcup_{B\in\Z_m(x_0,a)} B.
$$

\begin{proposition}\label{AmU}
Let $U,W$ be  open sets, $U\subset W\subset X$, $\nu\in \Mp$, $x_0\in\reald$,
and $a\in (0,1)$. For every 
$m\in\nat$, let $A_m$ be the {\rm(}finite{\rm)} union of all balls  
$B\in\Z_m(x_0,a)$ such that $B\subset U\cap B(0,m)$ and 
$\dist(B,\reald\setminus U)\ge 1/m$. Then
\begin{equation}\label{qx}
 \lim_{m\to\infty} 
\nu^{A_m\cup\wc}(q)
=\nu^{U\cup\wc}(q)\qquad\mbox{ for every }q\in\Px.
\end{equation} 
\end{proposition} 

For the proof we shall need the following lemma.

\begin{lemma}\label{rhat2}
Let $a\in (0,1)$, $x_0\in\reald$, $r>0$, and $\delta>0$.  
Then there exists $m_0\in\nat$ such that, for all $m\ge m_0$
 and   $x\in X$ satisfying $B(x,r)\subset X$, 
\begin{equation}\label{R2}
       R_1^{Z_m(x_0,a)\cap B(x,r)}(x)> 1-\delta.
         \end{equation} 
\end{lemma} 

\begin{proof} 

1. Let us first suppose that $d>\a$. For every subset $A$ of $\reald$,
 we define
$$
           u_A:=\inf\{u\colon u\ge 0\mbox{ superharmonic on }\reald,
                u\ge 1\mbox{ on }A\}
$$
(if $A$ is bounded, then $\hat u_A$ is the equilibrium potential
of $A$). Let $Z:=Z_1(0,a)$.
Obviously,
$$
  \inf u_{Z}(\reald)=\inf u_{Z}([0,1]^d),
$$
that is,  the continuous superharmonic function $u_{Z}$ admits
a minimum. Therefore $u_{Z}$ is constant. Since $u_{Z}=1$
on $Z$, we see that $u_{Z}$ is identically $1$. Consequently,
the sequence $(u_{Z\cap B(0,k)})_{k\in\nat}$ increases to $1$ 
locally uniformly on $\reald$ as  $k\uparrow \infty$. 
Given $\delta>0$, we hence may choose $k\in\nat$ such that 
\begin{equation}\label{R1k}
                 u_{Z\cap B(0,k)}> 1-\frac\delta 2 \on [0,1]^d.
\end{equation} 
There exists $K\ge k$ with $u_{Z\cap B(0,k)}<\delta/2$ on 
$\reald\setminus B(0,K)$. Then the function $v:=(u_{Z\cap B(0,k)}-\delta/2)^+$
 is subharmonic on
$\reald\setminus Z$ and vanishes outside $B(0,K)$.

Now let $m>(K+1)/\min\{r,1\}$.
There exists $z\in \frac 1m(x_0+\ganzd)$ with $x-z\in [0,\frac 1m]^d$. We define
$$
              w(y):=v(m(y-z))\qquad (y\in \reald).
$$
Then $w\le 1$, $w$ is subharmonic on $\reald\setminus Z_m(x_0,a)$, and $w=0$ on
$\reald\setminus B(z,K/m)$. Since $B(z,K/m)\subset B(x,r)\subset X$,
we see that the restriction of $w$ on $X$ vanishes outside $B(x,r)$ 
and hence $ w|_X\le R_1^{Z_m(x_0,a)\cap B(x,r)}$.
Since $m(x-z)\in [0,1]^d$, we know by (\ref{R1k}) that $w(x)\ge 1-\delta$.
Thus (\ref{R2}) holds.

2. Let us finally consider the remaining case $\a=d=2$ (classical case 
in the plane) and let $\sqrt 2/m\le r/2$.
There exists $z\in Z_m(x_0,a)$ such that $|x-z|<\sqrt 2/m$. Then 
$B(z,r/2)\subset B(x,r)$. We define
$$
 p_m(y):=\min\bigl\{1,\frac{\ln (r/2)-\ln|y-z|}
                 {\ln (r/2)+\ln m-\ln a}    \bigr\}\qquad(y\in B(x,r))
$$
Moreover, $p_m=1$ on $B(z,a/m)$, $p_m\le 0$ on $\partial B(x,r)$, and $p_m$
is harmonic on the open set $\overset{\,\,\circ}B(x,r)\setminus B(z,a/m)$.
Therefore  $$
R_1^{Z_m(x_0,a)\cap B(x,r)}\ge p_m \on B(x,r).
$$
Since 
$$
             \frac{\ln (r/2)-\ln|x-z|}{\ln (r/2)+\ln m-\ln a} 
        \ge   \frac{\ln (r/2)+\ln m - \ln\sqrt 2}{\ln (r/2)+\ln m-\ln a},
$$
we see that $p_m(x)>1-\delta$, if $m$ is sufficiently large.
 This finishes the proof.

\end{proof}

\begin{proof}[Proof of Proposition \ref{AmU}] Let 
 $q \in \Px$,  $\delta>0$, and $p:=q+\delta p_0$, where $p_0\in\Px$
such that $p_0>0$ and $\nu(p_0)\le 1$.
We choose an arbitrary
sequence $(K_n)$ of compact sets which 
is increasing to $U$. For the moment, let us fix $n\in\nat$.
There exists $0<r<\frac 12 \mbox{dist}(K_n,\reald\setminus U)$ such that, 
for every $x\in K_n$, $p>q(x)$ on $B(x,r)$. 
If $m\in\nat$ such that $m\ge 1/r$ and $K_n\subset B(0,m)$, then 
$Z_m(x_0,a)\cap B(x,r)$ is a subset of $ A_m$ and hence, for every $x\in K_n$, 
$$
      R_{p}^{A_m}(x)\ge 
 R_{q(x)}^{Z_m(x_0,a)\cap  B(x,r)}(x)=q(x) R_1^{Z_m(x_0,a)\cap  B(x,r)}(x).
$$
Using  Lemma \ref{rhat2} we hence obtain  $m_n\in\nat$
 such that, for every   $m\ge m_n$,
$$
      R_{p}^{A_m}> (1-\delta)q \qquad\mbox{ on }K_n.
$$
By the definition of reduced functions, this implies that,
for every $m\ge m_n$,
$$
 R_q^{U\cup \wc} +\delta p_0\ge R_p^{U\cup \wc} \ge
 R_{p}^{A_m\cup \wc}
\ge  (1-\delta)R_q^{K_n\cup\wc}
$$
and therefore
\begin{equation}\label{qp0}
 \rhat_q^{U\cup \wc} +\delta p_0 \ge \rhat_{p}^{A_m\cup \wc}
 \ge (1-\delta)\rhat_q^{K_n\cup\wc}.
\end{equation} 
If  $n\uparrow \infty$, 
then $\rhat_q^{K_n\cup\wc}\uparrow \rhat_q^{U\cup\wc}$  
whence
\begin{equation}\label{nq}
\nu^{K_n\cup\wc}(q)=\int \rhat_q^{K_n\cup\wc}\,d\nu \,
\uparrow \int R_q^{U\cup\wc}\,d\nu =\nu^{U\cup\wc}(q).
\end{equation} 
Since $\nu(p_0)\le 1$, (\ref{qp0}) and (\ref{nq}) imply that
(\ref{qx}) holds. 
\end{proof}

\section{Joint shrinking of disjoint small balls}

The following simple facts on iterated balayage will be used
again and again.
If $\nu\in\Mp$ and $A,\tilde A$ are closed sets such that 
$\tilde A\subset A\subset X$ and $\nu(A)=0$, 
then, 
\begin{equation}\label{iterated}
 \nu^{\tilde A}=
\nu^A|_{\tilde A}+ (\nu^A|_{A\setminus \tilde A})^{\tilde A}\ge
\nu^A|_{\tilde A}
\end{equation} 
(see \cite[VI.9.4]{BH} for a far more general statement). 

In terms of harmonic kernels $H_V$ for open subsets of $X$,
defined by $H_V(x,\cdot)=\vx^\vc$ for $x\in V$  and $H_V(x,\cdot)=\vx$
for $x\in V^c$, this can be expressed by
$$
H_UH_{\tilde U}=H_U,
$$
whenever $U,\tilde U$ are open subsets of $X$ satisfying $U\subset \tilde U$.
In the classical case this is equivalent to the following
property of the generalized Dirichlet solution, that is, the
Perron-Wiener-Brelot solution of the Dirichlet problem. 
If $f$ is a continuous  $\Px$-bounded function  on the boundary 
$\partial\tilde U$ of $\tilde U$ in $X$,
then the generalized Dirichlet solution $\tilde h$ for $\tilde U$ and $f$ 
coincides on $U$ with the generalized Dirichlet solution $h$ for~$U$ 
and the boundary function $g$, where $g=f$  on
 $\partial U\cap \partial \tilde U$ and $g=\tilde h$ on 
$\partial U\cap U$ (see \cite[Lemma 8.39]{helms}).

Moreover,  for every $\nu\in\Mp$ and for all closed sets $A,B$  in $X$,
 \begin{equation}\label{AB}
            (\nu^A)^B(B)\le \nu^B(B).
 \end{equation}  
Indeed, it suffices to notice that both measures $(\nu^A)^B$
and $\nu^B$ are supported by $B$ and that
$(\nu^A)^B(1)=\nu^A(\rhat_1^B)\le \nu(\rhat_1^B)=\nu^B(1)$ 
by (\ref{bal-est}).

We recall that, for  every $\eta\in [0,1]$ and every closed ball 
$B$ in $\reald$  having center~$x_B$ and radius $r_B$,
we  denote the ball obtained by shrinking $B$ with the
factor~$\eta$ by $B^\eta$, that is,
$$
            B^\eta:=x_B+\eta(B-x_B).
$$

\begin{lemma}\label{essential-1}
Let $A$ be the union of finitely many closed balls $B_1,\dots, B_m$
which are contained in $X$ and pairwise disjoint. For every $t\in
[0,1]^m$, let $A_t\subset A$ denote the union of the balls  $B_i^{t_i}$,
$1\le i\le m$. Moreover, let $\nu\in\Mp$ such that $\nu(A)=0$, 
$\gamma_1,\dots,\gamma_m\in \real_+$, and 
$$
 \Gamma:=\{t\in [0,1]^m\colon \nut(B_i)\le \gamma_i,\ 1\le
              i\le m\}.
$$
Then there exists  $s\in\Gamma$ such that $s\ge t$ for every 
$t \in \Gamma$. Moreover,  $\nus(B_i) =\gamma_i$ 
for every $i\in\{1,\dots,m\}$ such that $s_i<1$.
\end{lemma} 

\begin{proof} Let us note first that
  $\nut(B_i)=\nut(B_i^{t_i})$  for every $t\in\Gamma$
and for every $1\le i\le m$, since $\nut$ is supported by the
subset $A_t$ of $A$.

0. Of course, $(0,\dots,0)\in\Gamma$, since points are polar.

1. If $t,\tilde t\in\Gamma$, then $t\vee \tilde t\in \Gamma$. Indeed,
let us fix $1\le i\le m$. We may assume without loss of generality that 
$t_i\ge \tilde t_i$. Since $A_t\subset A_{t\vee \tilde t}$, 
we conclude by (\ref{iterated})
that
$$
              \nu^{A_{t\vee \tilde t}}(B_i^{t_i\vee \tilde t_i})
=             \nu^{A_{t\vee \tilde t}}(B_i^{t_i})
\le           \nut(B_i^{t_i})\le \gamma_i.
$$

2. For every $t<(1,\dots,1)$,
the set $A_t$ is the intersection of all $A_{\tilde t}$, $\tilde t>t$.
For every $t>(0,\dots,0)$ the set $A_t$ is the fine closure of the union of all
$A_{\tilde t}$, $\tilde t<t$. This implies that, for every $p\in \Px$, 
the mapping
$t\mapsto \nut(p)$ is continuous on $[0,1]^m$. Hence, for every 
$f\in\K(X)$, the mapping $t\mapsto \nut(f)$ is continuous.
Since the closed balls $B_1,\dots,B_m$ are disjoint, we obtain
that the  mapping 
$$
t\mapsto \left(\nut(B_1),\dots,\nut(B_m)\right)
$$
from $[0,1]^m$ into $[0,1]^m$ is continuous. Therefore $\Gamma$ is closed.

3. Combining (1) and (2) we see that 
$$
s:=(\sup_{t\in\Gamma} t_1,\dots,\sup_{t\in\Gamma} t_m)\in \Gamma,
$$
where of course $s\ge t$ for every $t\in\Gamma$.
To finish the proof, let us consider $i\in\{1,\dots,m\}$ such that $s_i<1$
and suppose that $\nus(B_i)<\gamma_i$. Let us define 
$\tilde s:=(s_1,\dots,s_{i-1},b,s_{i+1},\dots,s_m)$, where  $s_i<b\le 1$.
By (2),  we may choose $b$ in such a way that 
$\nu^{A_{\tilde s}}(B_i)<\gamma_i$. Since $A_s\subset A_{\tilde s}$,
we obtain by (\ref{iterated}) that 
$\nu^{A_{\tilde s}}(B_j^{s_j})\le \nus(B_j^{s_j})\le \gamma_j$
for every $j\in\{1,\dots,m\}$, $j\ne i$.
Thus $\tilde s\in \Gamma$, $\tilde s\le s$, $b=\tilde s_i\le s_i$,
a~contradiction.
\end{proof} 

Let us note the following simple consequence.

\begin{proposition}\label{Lempty}
Let $A$ be the union of finitely many closed balls $B_1,\dots, B_m$
which are contained in $X$ and pairwise disjoint. Moreover,  
let $\nu\in\Mp$ such that $\nu(A)=0$ and let $\b_1,\dots,\b_m\in [0,1]$. 
Then there exist $s_1,\dots,s_m\in[0,1]$
such that the union $\tilde A$ of the shrinked balls $B_1^{s_1},\dots,
B_m^{s_m}$ satisfies
$$
    \nu^{\tilde A}(B_i)=\b_i\nu^A(B_i) \qquad\mbox{ for every }1\le
    i\le m.
$$
\end{proposition} 

\begin{proof} It suffices to take $\gamma_i:=\b_i\nu^A(B_i)$, 
$1\le i\le m$, and to choose $s=(s_1,\dots,s_m)$ in $[0,1]^m$ 
according to Lemma
\ref{essential-1}. Then $\nu^{\tilde A}(B_i)\le \b_i\nu^A(B_i)$ for all $1\le
i\le m$. Furthermore, equality holds whenever $s_i<1$. If, however,
$i\in\{1,\dots,m\}$ such that 
$s_i=1$, then $\nu^{\tilde A}(B_i)\ge \nu^A(B_i)$ by (\ref{iterated})
whence as well $\nu^{\tilde A}(B_i)\ge \b_i \nu^A(B_i)$ (and $\b_i=1$
unless $\nu^A(B_i)=0$).
\end{proof} 
In the classical case $\a=2$, the harmonic measure $\vy^\uc$
for a ball $U=\overset\circ {\! B}(y_0,r)$ and $y\in U$ has the 
Poisson density
$$
        \rho_y^U(z)=r^{d-2}(r^2-|y-y_0|^2)|y-z|^{-d}, \qquad |z-y_0|=r,
$$
with respect to normalized surface measure on the boundary of $U$. 
For Riesz potentials (the case $0<\a<2$), $\vy^\uc$ has a density
$\rho_y^U$ with respect to Lebesgue measure on~$\uc$,
$$
        \rho_y^U(z)=a_\a\,\frac{(r^2-|y-y_0|^2)^{\a/2}}
                              {(|z-y_0|^2-r^2)^{\a/2}}\, |z-y|^{-d},
\qquad |y-y_0|<r\le |z-y_0|,
$$
where $a_\a$ is a constant depending on $d$ and $\a$
(see \cite[p.\,192]{BH}). 
If $y, \tilde y\in B(y_0,\eta r)$, $0<\eta<1$,
then in both cases 
\begin{equation}\label{poisson}
 \frac{\rho_y^U(z)}{\rho_{\tilde y}^U(z)}
\le \frac 1{(1-\eta^2)^{\a/2}}\frac {(1+\eta)^{d}}{(1-\eta)^{d}}
=\frac {(1+\eta)^{d-\frac \a 2}}{(1-\eta)^{d+\frac \a2}}.
\end{equation} 
If $\eta$ is small, then the expression on the right side of (\ref{poisson}) 
is approximately $1+2d\eta$. So there exists $\delta_0\in (0,1)$ such that 
\begin{equation}\label{harnack-null}
   \ve_{\tilde y}^\uc\le (1+\delta) \vy^\uc,
\end{equation} 
whenever $0<\delta\le \delta_0$ and $y,\tilde y\in
B(y_0,\frac{\delta}{3d}r)$. 
For the moment, let us fix  $\delta\in (0,\delta_0]$ and a closed subset
$A$ of $X$ such that $B(y_0,r) \cap A=\emptyset$. We observe that, 
for every  $y\in B(y_0,r)$, $\vy^A=(\vy^\uc)^A$. Indeed,
in the classical case this follows from  (\ref{iterated}), 
since then $\ve_y^\uc(A)=0$. In the general case,
it follows from $(\vy^A)^A=\vy^A$ (see \cite[VI.5.21]{BH}), since
trivially $(\vy^A)^A\prec (\vy^\uc)^A\prec \vy^A$. Therefore
(\ref{harnack-null}) implies that
\begin{equation}\label{harnack}
      \ve_{\tilde y}^A\le (1+\delta)\vy^A  \qquad
   \mbox{ for all }y,\tilde y\in B\bigl(y_0,\frac \delta{3d}r\bigr).
\end{equation} 
Let us say that finite family $\mathcal B$ of closed balls,
which are contained in $X$ and pairwise disjoint, is a
\emph{$\delta$-family in $X$}, if $0<\delta<\delta_0$ and the union $A$ of all
$B\in\mathcal B$ satisfies 
\begin{equation}\label{dist}
        r_B\le \frac\delta{3d}\,\mbox{dist}\bigl(x_{B},
        (\reald\setminus X)\cup
        (A\setminus B)\bigr) \quad\mbox{ for every }B\in\mathcal B.
\end{equation} 

Here is the  key to Theorem \ref{main-int}. As already indicated,
it will be applied to balayage on subsets $A$ of $W$
\emph{with respect to $W$} in place
of $X$ to deal with balayage measures of the
form $\nu^{A\cup\wc}$.

\begin{proposition}\label{essential-2}
Let $A$ be the union of a $\delta$-family  $B_1,\dots, B_m$ in $X$ and let
$\nu\in\Mp$ such that $\nu(A)=0$.
Moreover, let $\lambda\in\L$, let $I_1,\dots,I_k$ be a partition of 
$\{1,\dots,m\}$, 
and $K_n$ be the union of the balls $B_i$, $i\in I_n$, $1\le n\le k$.

Then there exist $s_1,\dots,s_m\in [0,1]$ such that 
$C:=B_1^{s_1}\cup\dots\cup B_m^{s_m}$ satisfies
\begin{equation}\label{alpha}
   \nu^C(B_i)= (1-\delta) \sulan \nun(B_i)
        \qquad \mbox{ for every } 1\le i\le m.\end{equation} 
\end{proposition} 

\begin{proof} Since the measures $\nun$ are supported by $K_n$,
the sum on the right side of (\ref{alpha}) reduces to the term $\lan\nun(B_i)$,
when $i\in I_n$. By Lemma \ref{essential-1}, there exists $s\in [0,1]^m$ 
such that $C:=B_1^{s_1}\cup\dots\cup B_m^{s_m}$ satisfies
\begin{equation}\label{KL-ine}
   \nu^C(B_i)\le  (1-\delta) \lan\nun(B_i)
        \qquad \mbox{ for all  } 1\le n\le k,\ i\in I_n,
\end{equation} 
with equality whenever $s_i<1$. We claim that we even have
\begin{equation}\label{contra}                 
 \nu^C(B_i)\ge  \lan\nun(B_i)\qquad\mbox{ if } s_i=1,\ i\in I_n,\ 1\le n\le k.
\end{equation} 
and this will clearly finish the proof (in fact, it shows even that
$s_i$ cannot be equal to $1$ for $i\in I_n$, unless 
$\lan\nun(B_i)=0$). 

Indeed, let us suppose, for example, that $s_l=1$
for some $l\in I_1$ and let $I_1':=I_1\setminus \{l\}$.
Then $B:=B_l=B_l^{s_l}$, that is, $B$ is a subset of $C$, and we get
by (\ref{iterated}) that
\begin{equation}\label{B-iterated}
\nu^B=\nu^C|_B+(\nu^C|_{C\setminus B})^B,
\end{equation} 
where
\begin{equation}\label{B-sum}
\nu^C|_{C\setminus B}=\sum_{i\in I_1'} \nu^C|_{B_i}+
                  \sum_{n=2}^k \sum_{i\in I_n} \nu^C|_{B_i}.
\end{equation} 
By (\ref{harnack}),
$$
   (\nu^C|_{B_i})^B\le (1+\delta)(1-\delta)\lambda_1
   (\nu^{K_1}|_{B_i})^B\le   \lambda_1 (\nu^{K_1}|_{B_i})^B
\qquad\mbox{ for all }i\in I_1'.
$$
Similarly, $(\nu^C|_{B_i})^B\le \lan (\nun|_{B_i})^B$ for all
$i\in I_n$, $2\le n\le k$. Taking sums  we see that 
$$
\sum_{i\in I_1'}  (\nu^C|_{B_i})^B \le \lambda_1 (\nu^{K_1}|_{K_1\setminus B})^B
\quad\mbox{ and}\quad
\sum_{i\in I_n} (\nu^C|_{B_i})^B
\le \lan (\nun)^B
$$
for every $2\le n\le k$.
Therefore (\ref{B-iterated}) and (\ref{B-sum}) imply the inequality
\begin{equation}\label{BKL}
\nu^B(B)\le \nu^C(B)+\lambda_1 (\nu^{K_1}|_{K_1\setminus B})^B(B)+
                    \sum_{n=2}^k\lan (\nun)^B(B),
\end{equation} 
where $(\nun)^B(B)\le \nu^B(B)$ by (\ref{AB}). Hence
$$
 \lambda_1\nu^B(B)\le \nu^{C}(B)+ \lambda_1 (\nu^{K_1}|_{K_1\setminus B})^B(B).
$$
Knowing that $  \nu^B=\nu^{K_1}|_B+(\nu^{K_1}|_{K_1\setminus B})^B$
by (\ref{iterated}), we thus get the inequality
$\lambda_1 \nu^{K_1}(B)\le \nu^C(B)$, and the proof is finished.
\end{proof} 

Let us note that a combination of Proposition \ref{Lempty} and
Proposition \ref{essential-2} yields the following stronger result
(which will not be needed in the sequel).

\begin{corollary}\label{essential-3} Under the assumptions of
  Proposition \ref{essential-2} and given any real numbers 
$\beta_1,\dots,\beta_m\in [0,1]$, there exist $\gamma_1,\dots,\gamma_m\in [0,1]$
such that the union $C$ of the shrinked balls
$B_1^{\gamma_1},\dots, B_m^{\gamma_m}$ satisfies
$$
   \nu^C(B_i)= (1-\delta) \beta_i\lambda_n\nun(B_i)
        \qquad \mbox{ for all } 1\le n\le k,\ i\in I_n.
$$
\end{corollary}

\section{Proof of Theorem \ref{main-int}}

1. Let $U_1,\dots,U_k$ be open subsets of an open set $W$ in $X$,
 let $\nu\in\Mp$ such that $\nu(\wc)=0$, and  $\lambda\in \L$.
We fix a strictly positive $p\in\Px$ such that $\nu(p)<\infty$.
To prove the weak convergence of a sequence 
$(\mu_n)$ to $\mu$, it is sufficient to check the convergence 
of the sequence 
$(\mu_n(f))$ to $\mu(f)$ for each functions $f$ from a suitable countable
subset of $\K(X)$. For every $f\in\K(X)$ and for every $\eta>0$, 
there exist $q,q'\in \Px$ which are bounded by a multiple of $p$
such that $|f-(q-q')|\le \eta p$ (see  \cite[I.1.2, III.6.10]{BH}).
To prove Theorem \ref{main-int}, it is therefore sufficient to show
the following. Let $\Q$ be a finite set of potentials $q\in\Px$ which 
are bounded by $p$ and let $0<\eta<1$.
Then there exists  a  union $C$ of pairwise disjoint closed balls in 
$U:=U_1\cup\dots \cup U_k$ such that 
\begin{equation}\label{suff-main}
        \bigl| \nu^{C\cup \wc}(q)- \sulan \nuun(q)\bigr| < \eta \qquad
\mbox{ for every }q\in Q.
\end{equation} 
By \cite[VI.1.9]{BH}, we may assume without loss of generality that 
all $\ov U_n$, $1\le n\le k$, are compact subsets of $W$ and that
$p\ge 1$ on $\ov U$. We  then define $\delta:=(6\nu(p)+3k)\inv\eta$.
 There exists $0<\delta'\le \delta$ such that 
\begin{equation}\label{continuous}
  |q(y)-q(z)|<\delta,\quad\mbox{ whenever }q\in\Q \mbox{ and }
                 y,z\in  \ov U,\ |y-z|<\delta'.
\end{equation}

2. By Proposition \ref{qx}, we are able to replace each $U_n$, $1\le n\le k$,
 by a~finite union~$K_n$ of very small closed balls. We then want to 
shrink  these balls  using Proposition~\ref{essential-2}. 
This, however, not only requires that $K_1,\dots,K_k$ be pairwise
disjoint, but that they are  separated well enough 
to obtain a $\delta$-family. 
Moreover, taking $A:=K_1\cup\dots\cup K_k$ we shall have to replace $\nu$
by the measure $\tilde \nu:=1_{X\setminus A}\nu$ not charging $A$,
and therefore $\nu(A)$ will have to be small. This can be achieved
 considering $\Z_m(x_0,a)$ for sufficiently many points $x_0\in\reald$.
We fix a natural number $N>k+\nu(p)/\delta$ and define 
$$
a:=\min\{\frac\delta{4dN}, \frac {\delta'}2\}, \qquad 
x_j:=\bigl(\frac jN,0,\dots,0\bigr), \quad 1\le j\le N.
$$
 For every  $M\in\nat$, $1\le n\le k$, and $1\le j\le N$, 
let $A_M(n,j)\subset Z_M(x_j,a) $ be the union of all balls 
$B\in\Z_M(x_j,a)$ such that 
$B\subset  U_n$ and $\dist(B,\reald\setminus  U_n)\ge 1/M$.
By Proposition \ref{AmU}, there exists $M\in\nat$ such that,
for every $q\in\Q$ and for all $1\le n\le k$, $j\in\{1,\dots,N\}$,
\begin{equation}\label{estW}
  |\nu^{A_M(n,j)\cup\wc}(q)-\nu^{U_n\cup\wc}(q)| < \delta .
\end{equation} 
By our definition of $a$ and $x_1,\dots,x_N$, the sets
$Z_M(x_1,a),\dots, Z_M(x_N,a)$ are pairwise disjoint and hence
$$
     \sum_{j=1}^N\nu(1_{Z_M(x_j,a)}p)\le \nu(p)<(N-k)\delta.
$$
Therefore at least $k$ of terms of the sum must be strictly smaller than
$\delta$, that is, there exist $k$ different 
$j_1,\dots,j_k\in\{1,\dots,N\}$ such that 
\begin{equation}\label{nu-tilde-nu}
      \nu(1_{Z_M(x_{j_n},a)}p)<\delta \qquad\mbox{ for every } 1\le n\le k.
\end{equation} 
We define 
$$
   K_1:=A_M(1,j_1),\ \dots\ ,   K_k:=A_M(k,j_k),\quad 
   A:=K_1\cup\dots\cup K_k, \quad\tilde\nu=1_{X\setminus A}\nu.
$$
 By (\ref{nu-tilde-nu}), $(\nu-\tilde\nu)(p)<k\delta$. 
The set $A$ is a union of pairwise disjoint balls $B_1,\dots,B_m$ from
the union of  $\Z(x_{j_1},a)$, $\dots$, $\Z(x_{j_k},a)$. 
Hence, for every $1\le i\le m$,
$$   \dist\bigr(x_{B_i},(\reald\setminus W)\cup (A\setminus B_i)\bigr)
\ge \frac 1M\bigl(\frac 1N-{2a}\bigr)\ge\bigl(\frac{4d}\delta-2\bigr)\frac aM
\ge \frac{3d}\delta\frac aM= \frac{3d}\delta r_{B_i}.
$$
So $B_1,\dots,B_m$ is a $\delta$-family in $W$ and, of course, 
 $\tilde \nu(A)=0$. Thus we may apply Proposition \ref{essential-2}
to~$W$ in place of $X$ and to $\tilde \nu$ in place of $\nu$.
Denoting balayage of $\tilde \nu$ on compact  subsets $L$ of $W$ 
\emph{relative to $W$}
by $\wtnu^L$, we obtain $s_1,\dots, s_m \in [0,1]$ such that the union
$C$ of the shrinked balls $B_1^{s_1},\dots ,B_m^{s_m}\subset A$ satisfies 
$$
          \wtnu^{C}(B_i)= (1-\delta) \sulan \wtnu^{K_n} (B_i)\qquad
   \mbox{ for every }1\le i\le m.
$$
By \cite[VI.2.9]{BH}, this means that defining
$\mut:=\sulan\nutkn$
we have
\begin{equation}\label{alpha-proposition}
   \tilde\nu^{C\cup \wc}(B_i)=(1-\delta) \mut(B_i) \qquad
\mbox{ for every }1\le i\le m.
\end{equation} 
 
3. We now fix $q\in \Q$ and claim first that
\begin{equation}\label{CKL}
 \bigl|\tilde \nu^{C\cup \wc}(1_W q)-\mut(1_W q)\bigr|
< 4\nu(p)\delta.
\end{equation} 
Indeed, let $ g:=\sum_{i=1}^m q(x_{B_i}) 1_{B_i}$. 
By (\ref{alpha-proposition}),
\begin{equation}\label{CKLg}
\tilde\nu^{C\cup\wc}(g)= (1-\delta)\mut(g).
\end{equation} 
Since $0\le g\le 2p$, we know by (\ref{bal-est}) that, for every 
$1\le n\le k$,
$$
0\le \nutkn(g)\le2 \nutkn(p) \le 2\nu(p).
$$
Moreover, $|g-1_W q|< \delta p$ on $\ov{U}\cup \wc$ by  (\ref{continuous}).
Therefore
$$
\bigl|\tilde\nu^{C\cup \wc}(g)-\tilde\nu^{C\cup\wc}(1_W q)\bigr|
\le \delta\tilde\nu^{C\cup\wc}(p)\le \nu(p) \delta 
$$
and, for every $1\le n\le k$,
$$
\bigl|\nutkn(g)-\nutkn(1_W q)\bigr|<\delta\nutkn(p) \le \nu(p)\delta.
$$
Thus (\ref{CKLg}) implies (\ref{CKL}) (and the proof would readily
 be  finished in the case $W=X$).

4. It may be surprising that  (\ref{alpha-proposition}), which merely
indicates 
that $\tilde \mu$ is a good approximation for $\tilde\nu^{C\cup\wc}$ on $W$, 
also implies
that $\tilde \mu$ approximates $\tilde \nu^{C\cup\wc}$ nicely on~$X\setminus W$.
We claim that 
\begin{equation}\label{tau}
\rho:=\tilde \nu^{C\cup\wc}|_\wc-\mut|_\wc \ge 0, \quad 
\mbox{ and }\quad \rho(p)\le 2\nu(p)\delta.
\end{equation} 
Indeed, by (\ref{iterated}), 
$$
 \tilde \nu^{C\cup \wc}|_\wc+(\tilde\nu^{C\cup \wc}|_W )^\wc
=\tilde\nu^\wc=  \mut|_\wc+(\mut|_W)^\wc.
$$
Defining $\sigma:=\mut|_W$ and $ \tau:=\tilde\nu^{C\cup\wc}|_W$
we hence see that 
$$
\rho=\sigma^\wc -\tau^\wc.
$$
By (\ref{alpha-proposition}) and (\ref{harnack}), 
for each $B\in\{B_1,\dots,B_m\}$,
$$
(1_B\tau)^\wc\le (1-\delta)(1+\delta)(1_B\sigma)^\wc
            \le  (1_B\sigma)^\wc,
$$
$$
(1_B\tau)^\wc\ge \frac {1-\delta}{1+\delta}  (1_B\sigma)^\wc
\ge (1-2\delta) (1_B\sigma)^\wc.
$$
Taking the sum we obtain that
$
      0\le \rho\le  2\delta  \sigma^\wc
$
where $\sigma^\wc(p)\le \mut^\wc(p)\le \mut(p)\le\nu(p)$ by (\ref{bal-est}).
Thus (\ref{tau}) holds and we conclude that
\begin{equation}\label{CKLc}
 \bigl|\tilde \nu^{C\cup \wc}(1_\wc q)-
\mut (1_\wc q)\bigr|=\rho(q)\le\rho(p)
\le 2\nu(p)\delta.
\end{equation} 

5. Combining (\ref{CKL}) and (\ref{CKLc}) we get
$
 \bigl|\tilde \nu^{C\cup \wc}(q)-\mut ( q)\bigr|< 6 \nu(p)\delta
$.
This implies 
$$
 \bigl| \nu^{C\cup \wc}(q) - \sulan \nukn( q)\bigr|< 6\nu(p)+2k\delta,
$$
since $\tilde \nu\le \nu$ and $(\nu-\tilde \nu)(p)<k\delta$.
Together with (\ref{estW}), this estimate finally yields 
$$
 \bigl| \nu^{C\cup\wc}(q) -  \sulan\nuun (q)\bigr|
 <6\nu(p)\delta+3k\delta= \eta,
$$
that is, (\ref{suff-main}) holds and the proof is finished.
 
\section{Proof of Theorem \ref{jensen-char}}

Let us return to classical potential theory. We fix an
open subset $\Omega$ of $\reald$, $d\ge 2$, and a point $x\in\Omega$.
If $W$ is a bounded open subset of $\Omega$, then the measures in
$\mxs$ are supported by $\ov W$ and $\mxs$ is independent 
of the choice of the Greenian domain $X$ containing $W$
(of course, we would simply take $X=\Omega$ if $d\ge 3$ or if, in the case
$d=2$, the complement of $\Omega$ is non-polar).

\begin{proposition}\label{jx}
$$
  J_x(\Omega)=\bigcup\bigl\{ \mxs\colon W\mbox{ open},\ x\in W,\
                             \ov W\mbox{ compact in } \Omega\bigr\}.
$$
\end{proposition} 

\begin{proof} 1. Let $W$ be a bounded open set such that $x\in W$,
 $\ov W \subset \Omega$, and let $\mu\in\mxs$. Then (in
contrast to the situation for Riesz potentials) $\mu$ is supported by
$\ov W$. Let $X$ be a bounded domain such that $\ov W\subset X\subset
\Omega$. If $v$ is a continuous superharmonic function on $\Omega$,
then there exists a $\Px$-bounded continuous function on $X$ such that
$\tilde v=v$ on $\ov W$ whence $\tilde v\in S(W)$ and 
$
       \int v\,d\mu=\int \tilde v\,d\mu\le \tilde v(x)= v(x)
$.
Thus $\mu\in J_x(\Omega)$.

2. Let us now suppose conversely that $\mu\in J_x(\Omega)$. Let $K$
be a compact neighborhood of the  support
of $\mu$ such that $x\in K$ and each bounded component of $\reald\setminus K$
meets $\reald\setminus \Omega$. Further, let $W$ be a bounded open
neighborhood of $K$ such that $\ov W\subset \Omega$ and $v\in S(W)$.
By \cite[Theorem 6.1]{gardiner-app}, there exists a continuous
superharmonic function~$\tilde v$ on~$\Omega$ such that $\tilde v=v$
on $K$ whence
$
     \int v\,d\mu=\int \tilde v\,d\mu\le \tilde v(x)=v(x)
$.
Thus $\mu\in\mxs$.
\end{proof} 

A consequence is a full characterization of extremal Jensen
measures which has been asked for in \cite{cole-ransford}.

\begin{corollary}\label{} The set of all extremal Jensen measures for
$x$ with respect to~$\Omega$ is given by
\begin{eqnarray*}
  \bigl(J_x(\Omega)\bigr)_e & =&
\bigcup\bigl\{ \mxse\colon W\mbox{ open},\ x\in W,\
                             \ov W\mbox{ compact in } \Omega\bigr\}\\[2mm]
&=&\bigl\{\vx^{A^c}\colon \ov A \mbox{ compact in }\Omega\bigr\}.
\end{eqnarray*} 
\end{corollary} 

\begin{proof} 
1. If $W\subset \tilde W$, then $S(\tilde W)\subset S(W)$,
$\mxs\subset \M_x(S(\tilde W))$, and $\mxs$ is a closed face
of the compact convex set $\M_x(S(\tilde W))$ (see Theorem~\ref{msg}).
By Proposition  \ref{jx}, this immediately yields the first identity.

2. To prove the second identity, let $\mu\in\mxse$ for some open 
set $W$ such that $x\in W$ and
$\ov W$ is a compact subset of $\Omega$.
By Theorem \ref{msg}, $\mu=\vx^B$ for some set $B$ containing the
complement of $W$. Taking $A:=B^c$ we have $\mu=\vx^{A^c}$ and
$\ov A\subset \ov W\subset \Omega$.

 Conversely, let $A$ be a bounded set in $\reald$ such that $\ov
A\subset \Omega$ and consider $\mu=\vx^{A^c}$. 
Let $W$ be a bounded open set such that $x\in W$,
$\ov A\subset W$, and $\ov W\subset \Omega$. Then $\mu\in \mxse$ by
Theorem \ref{msg}.
\end{proof}

\section{Appendix}

Finally, let us give a proof for Theorem \ref{msg} in the 
general case. 
For the moment, we fix a closed set $A$ in  $X$ and recall
that the base $b(A)$ of  $A$
is the set of all points  $x\in A$ such that $A$ is not thin at $x$,
that is, $\vx^A= \vx$. 
Since in our case of classical potential theory or Riesz potentials
every semi-polar set is polar, the set $A\setminus b(A)$ is polar, hence
$b(b(A))=b(A)$ and $\mu^{b(A)}=\mu^A$ for every $\mu\in\Mp$
(see \cite[VI.5.12]{BH}). Moreover, $\beta(A)$, which can be
characterized as being the largest subset~$\tilde A$ of 
$b(A)$ such that $b(\tilde A)\subset \tilde A$, coincides with $b(A)$ (see
\cite[VI.6.1, VI.6.6]{BH}). Hence $\mu^{\beta(\wc)}=\mu^\wc$ for every 
$\mu\in\Mp$ and every open $W$ in $X$.

\begin{theorem}\label{msg-app}
Let $W$ be an open subset of $X$ and  $\nu\in\Mp$
such that $\nu(\wc)=0$. Then
\begin{eqnarray*}
\mns&=&\mnp\cap\mnh\\&=&\{\mu\in\mnp\colon \mu^\wc=\nu^\wc\}
                  =\{\mu\in\Mp\colon \nu^\wc\prec\mu\prec \nu\}.
\end{eqnarray*} 
Moreover, $\mns$ is a closed face of $\mnp$ and
$$
     \mnse= \{\onu^A\colon \wc\subset A\subset X,\ A \mbox{ Borel set}\}.
$$
\end{theorem} 

\begin{proof} 
Replacing the measure $\vx$ in the proof
of \cite[VI.9.5]{BH}  by $\nu$ (and  using that
$\mu^{\beta(\wc)}=\mu^{\wc}$),  we obtain immediately the first 
three identities, the fact that $\mns$ is a closed face of $\mnp$, and that
every measure $\onu^A$, where $A$ is a Borel set with 
$\wc\subset A\subset X$, is contained in $\mnse$.

Conversely, let $\mu\in\mnse$. Of course, 
$\mu\in\mne$, since the set $\mns$
is a closed face of $\mnp$. So, by (\ref{ext}), there exists a Borel 
subset~$A$ of~$X$ such that $\mu=\onu^A$. 
We intend to show that $\mu=\onu^{A\cup\wc}$. This will finish the proof,
since $A\cup \wc$ is a Borel subset of $X$ containing $\wc$.

By the characterization of $\mns$ given above, $\nu^\wc \prec\mu$. 
By \cite[VI.1.9]{BH}, this implies that,  for every $p\in \Px$,
$$
 \nu^\wc(p)=\nu^\wc( R_p^\wc)=
\inf_{U\mbox{\tiny open }\supset \wc}\nu^\wc(R_p^U)\le
\inf_{U\mbox{\tiny open }\supset \wc}\onu^A( R_p^U)=
\onu^A(R_p^\wc),
$$
 that is, 
\begin{equation}\label{final1}
 \nu^\wc\prec \onu^A|_\wc+(\onu^A|_W)^\wc
\end{equation} 
(see the proof of \cite[VI.9.9]{BH}). In addition,
\begin{equation}\label{final2}
\onu^{A\cup \wc}+\onu^A|_\wc+(\onu^A|_W)^\wc\prec \onu^A+\nu^\wc.
\end{equation} 
Indeed, if $\nu(A)=0$, this follows from \cite[VI.9.8]{BH}.
And if $\nu(A^c)=0$, then $\onu^{A\cup\wc}=\onu^A=\nu$ and 
(\ref{final1}) reduces to the trivial
statement $\nu+\nu^\wc\prec\nu+\nu^\wc$. The general case follows
decomposing $\nu$ into $1_{A^c}\nu$ and $1_A\nu$.

Combining (\ref{final1}) and (\ref{final2}), we see that $\onu^{A\cup
  \wc}\prec \onu^A$. Since trivially $\onu^A\prec \onu^{A\cup \wc}$,
we conclude that $\mu=\onu^A= \onu^{A\cup \wc}$ as claimed above, and
the proof is finished.
\end{proof}

  \bibliography{lit_bank}

\begin{thebibliography}{10}

\bibitem{baxter-chacon}
J.~R. Baxter and R.~V. Chacon.
\newblock Compactness of stopping times.
\newblock {\em Z. Wahrscheinlichkeitstheorie und Verw. Gebiete},
  40(3):169--181, 1977.

\bibitem{BH}
J.~Bliedtner and W.~Hansen.
\newblock {\em {Potential Theory -- An Analytic and Probabilistic Approach to
  Balayage}}.
\newblock Universitext. Springer, Berlin-Heidelberg-New York-Tokyo, 1986.

\bibitem{cole-ransford}
B.~J. Cole and T.~J. Ransford.
\newblock Jensen measures and harmonic measures.
\newblock {\em J. Reine Angew. Math.}, 541:29--53, 2001.

\bibitem{falkner}
N.~Falkner.
\newblock Stopped distributions for {M}arkov processes in duality.
\newblock {\em Z. Wahrsch. Verw. Gebiete}, 62(1):43--51, 1983.

\bibitem{falkner-fitzsimmons}
N.~Falkner and P.~J. Fitzsimmons.
\newblock Stopping distributions for right processes.
\newblock {\em Probab. Theory Related Fields}, 89(3):301--318, 1991.

\bibitem{fitzsimmons-sko}
P.~J. Fitzsimmons.
\newblock Skorokhod embedding by randomized hitting times.
\newblock In {\em Seminar on Stochastic Processes, 1990 (Vancouver, BC, 1990)},
  volume~24 of {\em Progr. Probab.}, pages 183--191. Birkh\"auser Boston,
  Boston, MA, 1991.

\bibitem{gardiner-app}
S.~J. Gardiner.
\newblock {\em Harmonic approximation}, volume 221 of {\em London Mathematical
  Society Lecture Note Series}.
\newblock Cambridge University Press, Cambridge, 1995.

\bibitem{helms}
L.~L. Helms.
\newblock {\em Introduction to potential theory}.
\newblock Pure and Applied Mathematics, Vol. XXII. Wiley-Interscience, New
  York-London-Sydney, 1969.

\bibitem{moko-ext}
G.~Mokobodzki.
\newblock \'{E}l\'ements extr\'emaux pour le balayage.
\newblock In {\em S\'eminaire de Th\'eorie du Potentiel, dirig\'e par M.
  Brelot, G. Choquet et J. Deny (1969/70), Exp. 5}, page~14. Secr\'etariat
  Math., Paris, 1971.

\bibitem{ransford}
T.~Ransford.
\newblock {\em Potential theory in the complex plane}, volume~28 of {\em London
  Mathematical Society Student Texts}.
\newblock Cambridge University Press, Cambridge, 1995.

\bibitem{rost}
H.~Rost.
\newblock The stopping distributions of a {M}arkov {P}rocess.
\newblock {\em Invent. Math.}, 14:1--16, 1971.

\end{thebibliography}
\bibliographystyle{plain}

\vskip1cm

\centerline{
\vtop{\hsize=8.5truecm
\noindent Wolfhard Hansen\newline
Fakult\"at f\"ur Mathematik\newline
Universit\"at Bielefeld\newline
33501 Bielefeld\newline
Germany\\[2mm]
 hansen$@$math.uni-bielefeld.de}
\vtop{\hsize=7truecm
\noindent Ivan Netuka\newline
Charles University\newline
Faculty of Mathematics and Physics\newline
Mathematical Institute\newline
 Sokolovsk\'a 86\newline
 186 75 Praha 8\newline
 Czech Republic\\[2mm]
netuka@karlin.mff.cuni.cz}
\hfill
            }                                      
\end{document}